 \documentclass[journal]{IEEEtran}

\ifCLASSINFOpdf
\else
\fi

\usepackage{svg}
\usepackage{xcolor}
\usepackage{amssymb}
\usepackage[mode=buildnew]{standalone}
\usepackage{enumerate}
\usepackage{amsmath}
\usepackage{mathtools}
\usepackage{mathrsfs}
\usepackage{tabularx}
\usepackage{comment}
\usepackage{graphicx}
\usepackage{psfrag}
\usepackage{epsfig}
\usepackage{multicol}
\usepackage{multirow}
\usepackage{graphics}
\usepackage{subcaption}
\usepackage{algorithm}
\usepackage{algcompatible}
\usepackage{booktabs}
\usepackage{float}
\usepackage{arydshln}
\usepackage{multirow}
\usepackage{multicol}
\usepackage{graphicx}
\usepackage[hidelinks]{hyperref}
\usepackage{cite}
\usepackage{lipsum}
\usepackage{optidef}

\usepackage{color}

\newtheorem{theorem}{Theorem}

\usepackage{physics}

\begin{document}
\title{Gradient-Based Stochastic Extremum-Seeking Control for Multivariable Systems with Distinct Input Delays}    
	
	\author{Paulo Cesar Souza Silva, \IEEEmembership{Graduate Student Member, IEEE}, Paulo César Pellanda, \IEEEmembership{Member, IEEE}, \newline and Tiago Roux Oliveira, \IEEEmembership{Senior Member, IEEE.} 
		\thanks{P. C. S. Silva and P. C. Pellanda are with the Department of Defense Engineering, Military Institute of Engineering (IME), Rio de Janeiro,
		RJ 22290-270, Brazil (e-mails: cesar.paulo151@hotmail.com;  pcpellanda@ieee.org).}
     \thanks{T. R. Oliveira is with the Department of Electronics and
Telecommunication Engineering, State University of Rio de Janeiro (UERJ), Rio de Janeiro, RJ 20550-900, Brazil (e-mail: tiagoroux@uerj.br).}
		 \thanks{This work was supported in part by the Brazilian funding agencies CAPES---finance code 001, CNPq, and FAPERJ.}}

	\markboth{}%
	{Shell \MakeLowercase{\textit{et al.}}: Bare Demo of IEEEtran.cls for IEEE Journals}

	\maketitle

	\begin{abstract}
            This paper deals with the gradient extremum seeking control for static scalar maps with actuators governed by distributed diffusion partial differential equations (PDEs). To achieve the real-time optimization objective, we design a compensation controller for the distributed diffusion PDE via backstepping transformation in infinite dimensions. A further contribution of this paper is the appropriate motion planning design of the so-called probing (or perturbation) signal, which is more involved than in the non-distributed counterpart. Hence, with these two design ingredients, we provide an averaging-based methodology that can be implemented using the gradient and Hessian estimates. Local exponential stability for the closed-loop equilibrium of the average error dynamics is guaranteed through a Lyapunov-based analysis. By employing the averaging theory for infinite-dimensional systems, we prove that the trajectory converges to a small neighborhood surrounding the optimal point. The effectiveness of the proposed extremum seeking controller for distributed diffusion PDEs in cascade of nonlinear maps to be optimized is illustrated by means of numerical simulations.
	\end{abstract}


	\IEEEpeerreviewmaketitle

\section{Introduction}
\sloppy

\label{sec:introduction}
Extremum Seeking Control (ESC) is a model-independent, real-time adaptive optimization method aimed at determining the extremum point (maximum or minimum) of a nonlinear map \cite{KW:00}. Recently, significant advancements have been made in both the theory and application of ESC, including proofs of local \cite{KW:00,AK:03,K:2014} and semi-global \cite{TNM:06} stability of the search algorithm, even in the presence of local extrema \cite{TNMA:09}. Other developments include its extension to multivariable cases \cite{GKN:2012}, as well as improvements in parameter convergence and performance \cite{GZ:2003,AG:2007,GD:2015,NMM:2013}. The book \cite{LK:2012} further introduces stochastic versions of the algorithm, incorporating filtered noise perturbation signals.

On the other hand, actuator and sensor delays are common challenges in engineering practice. In this context, the backstepping transformation in infinite dimensions \cite{c1} has introduced a new approach to designing predictor feedback for delay compensation, allowing for the derivation of explicit Lyapunov functions for stability analysis in predictor-based control schemes.

In \cite{c8}, the authors pioneered a solution to the problem of multivariable \textit{deterministic} ESC algorithms for systems with output and/or equal input delays, using predictor feedback. They presented two approaches for implementing a predictor that compensates for delays with perturbation-based (or averaging-based) estimates of the model. The first approach employs gradient optimization, estimating the gradient and Hessian -- the first and second derivatives \cite{GKN:2012, c12} -- of the map to be optimized. The second approach is based on Newton optimization, where the Hessian’s inverse is estimated to make the convergence rate independent of the map’s unknown parameters. This setup presents an additional challenge, as predictor feedback typically requires a known model, whereas in these optimization problems, the maps are assumed to be unknown.

In \cite{c7}, the results are extended to systems with multiple and distinct input delays, applying both the gradient and Newton algorithms. This extension, particularly for the gradient case, leads to a complex predictor design. Stability analysis is conducted using successive backstepping transformations and the theory of averaging in infinite dimensions \cite{c9, c10}. Alternatively, references \cite{A16} and \cite{Emilia_arxiv} propose new designs for deterministic multivariable extremum seeking with arbitrarily long time delays, using finite-dimensional sequential predictors and sampled-data designs, respectively. These solutions eliminate the need for distributed integral terms but may result in delay-dependent conditions or affect the algorithm’s convergence rate, as ``larger delays or higher-dimension maps require smaller values of the dither excitation frequencies and therefore lead to slower convergence''.

In this paper and its companion (deterministic) conference version \cite{c21}, we develop an alternative multivariable gradient-based \textit{stochastic} ESC design for systems with \textit{multiple} and \textit{distinct} input delays, without the need for backstepping transformations. Unlike the approach in \cite{c7}, this design yields a significantly simpler feedback structure. This work represents the first predictor-feedback extremum-seeking result (with distributed-integral terms) that employs the classical reduction approach \cite{artstein} instead of backstepping transformations \cite{c1}. Furthermore, a key advantage over the results in \cite{A16} and \cite{Emilia_arxiv} is that our method achieves delay independence with respect to the dither frequency and system dimension. This means that achieving faster convergence rates does not require restricting the delay length or reducing the system’s order. In our method, there are no such restrictions on delay length, which can indeed be arbitrarily long. This aspect is crucial, as it preserves one of the primary advantages of stochastic extremum controllers over their deterministic counterparts: the potential for faster convergence rates \cite{LK:2012}. Our design aims to retain this benefit even when large time delays are present in the closed-loop system.

\begin{figure*}[htb!]
\centering \includegraphics[scale=0.75]{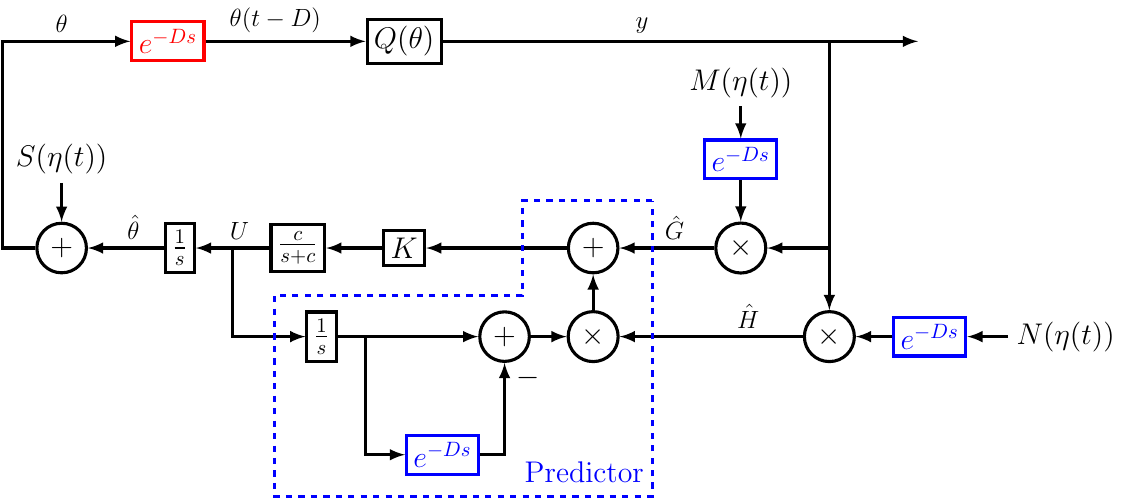}
\caption{Block diagram of the prediction scheme for compensating multiple and distinct input delays in the map $Q(\theta)$. The vector signals $\hat{G}$ and $\hat{H}$ represent the estimates for the gradient and Hessian of $Q(\theta)$, respectively. The multiple delays and gains are compactly represented by $D=$ diag$\{D_1,D_2,...,D_n\}$ and $K=$ diag$\{K_1,K_2,...,K_n\}$. 
The red block indicates the introduced delays, while the blue blocks show modifications to the classical stochastic gradient-based ESC algorithm \cite{LK:2012,c16} for mitigating the effects of time delays. In particular, the prediction feedback law is governed by (\ref{eq:14}), while the stochastic perturbations $S(t)$ and $M(t)$ and the demodulation signal $N(t)$ are given by (\ref{eq:4}), (\ref{eq:5}), and (\ref{eq:9}), respectively.}
\label{Fig1}
\end{figure*}


\section{Preliminaries}
 
 The partial derivatives of $u(x,t)$ are denoted by $u_t(x,t)$ and $u_x(x,t)$ or, when referring to the averaging operation $u_{av}(x,t)$, by $\partial_tu_{av}(x,t)$ and $\partial_xu_{av}(x,t)$.

 Consider a generic nonlinear system $\dot{x}=f(t,x,\epsilon)$, where the state $x \in$ $\mathbb{R}^n$, $f(t,x,\epsilon)$ is periodic with period $T$, such that $f(t+T,x,\epsilon)=f(t,x,\epsilon)$. For sufficiently small $\epsilon > 0$, the averaged model can be obtained as $\dot{x}_{av}=f_{av}(x_{av})$, with $f_{av}=\frac{1}{T}\int_{0}^{T} f(\tau,x_{av},0), d\tau$, where $x_{av}(t)$ represents the averaged version of the state $x(t)$ \cite{c13}.

 The 2-norm of a finite-dimensional state vector $X(t)$ in a system governed by a ordinary differential equation (ODE) is denoted by single bars, $|X(t)|$. In contrast, norms of functions (of $x$) are denoted by double bars. We denote the spatial $\mathcal{L}_2[0,D]$ norm of the PDE state $u(x,t)$ as $\|u(t)\|_{\mathcal{L}_2([0,D])}^2 := \int_{0}^{D}u^2(x,t)dx$. For simplicity, we omit the subscript $\mathcal{L}_2([0,D])$ in the following sections, using $\|\cdot\|$ in place of $\|\cdot\|_{\mathcal{L}_2([0,D])}$ unless otherwise specified \cite{c1}. The acronym PDE refers to partial differential equations.
 
 As defined in \cite{c13}, a vector-valued function $f(t,\epsilon) \in \mathbb{R}^n$ is said to be of order $O(\epsilon)$ within the interval $[t_1,t_2]$ if there exist positive constants $k$ and $\epsilon^{\ast}$, such that $|f(t,\epsilon)| \leq k \epsilon$,  $\forall \epsilon \in [0,\epsilon^\ast]$ and $\forall t \in [t_1, t_2]$.
 
 In this paper, $e_i \in \mathbb{R}^n$ denotes the $i^{th}$ column of the identity matrix $I_n \in \mathbb{R}^{n \times n}$.
 
Given an $\mathbb{R}^n$-valued signal $f$, the notation $f^D$ denotes
\begin{equation} \label{eq:2}
    f^D(t)=\big[f_1(t-D_1)\; f_2(t-D_2)\; ...\; f_n(t-D_n)\big]^T.
\end{equation}
Equation (\ref{eq:2}) is defined by the computation $f^D(t):=f(t-D)=e^{-Ds}[f(t)]$, where each input has a distinct, known, and constant delay.
These delays are ordered, without loss of generality, as
\begin{equation} \label{eq:1}
    D=\textrm{diag}\{D_1, D_2, \cdot \cdot \cdot, D_n \}, \;\; 0 \leq D_1 \leq D_2 \leq \ldots \leq D_n,
\end{equation}
with $s$ representing the Laplace operator.


\section{A Simplified Predictor for Gradient-Based ESC with Multiple and Distinct Input Delays}

Let $Q: \mathbb{R}^n \rightarrow \mathbb{R}$ be a convex static map with a maximum at $\theta^\ast \in \mathbb{R}^n$. It is assumed that the optimal input parameter $\theta^\ast$ is unknown, while the output $Q$ is available past input signals. Specifically, the measured signal is: 
\begin{equation} \label{eq:3}
    y(t)=Q(\theta^D(t)),
\end{equation}
\noindent
where 
\begin{equation} \label{TAChoje1}
    Q(\theta)=y^\ast+\frac{1}{2}(\theta-\theta^\ast)^TH(\theta-\theta^\ast),
\end{equation}
\noindent
with $y^\ast \in \mathbb{R}$ denoting the extremum point and $H=H^T<0$ representing the unknown $n\times n$ Hessian matrix of the static map.

The goal of extremum seeking is to estimate the parameter $\theta^\ast$ from the output $y$. From Figure~\ref{Fig1}, we have
\begin{equation} \label{eq:6}
    \hat{G}(t)=M(\eta^D(t))y(t)
    ~~~~\text{and}~~~~ \theta(t)=\hat{\theta}(t)+S(\eta(t)),
\end{equation}
where $\hat{\theta}$ denotes the estimate of $\theta^\ast$. 

The sinusoidal perturbation signals, referred to as dither signals, $S(\eta(t))$ and $M(\eta(t))$ $\in \mathbb{R}^n$, are defined as follows: 
\begin{equation} \label{eq:4}
    S(\eta(t))=\bigg[a_1\text{sin}(\eta_1(t)) \: ... \: a_n\text{sin}(\eta_n (t)) \bigg]^T\;, 
\end{equation}
\begin{equation} \label{eq:5}
    M(\eta(t))=\bigg[\frac{2}{a_1}\text{sin}(\eta_1 (t)) \: ... \: \frac{2}{a_n}\text{sin}(\eta_n (t))\bigg]^T\;,
\end{equation}
\noindent
where $a_i$, $i=1,2,...,n$, are nonzero amplitudes.

The elements of the stochastic Gaussian perturbation vector $\eta_i(t)$ are sequentially and mutually independent, such that $\mathbb{E}\{\eta(t)\}\!=\!0$,
$\mathbb{E}\{\eta^2_i(t)\}=\sigma^2_i$ and $\mathbb{E}\{\eta_i(t)\;\eta_j(t)\}=0$, $\forall i \neq j$, where $\mathbb{E}\{\cdot\}$ denotes the expectation of a signal. Additionally, we assume that the probability density function of the perturbation vector is symmetric about its mean. 
We apply a stochastic perturbation to the unknown map by incorporating the sinusoid of a Wiener process $W_{\omega t}$ around the boundary of a circle \cite{LK:2012,c14}. 
To ensure the Markov property, predictions of future states depend only on the current state of the process, independent of past states. Accordingly, we refer to this as a Markov process -- a stochastic process that satisfies the Markov property with respect to its natural filtration. 
%

The estimate of the unknown Hessian $H$ is given by 
\begin{equation} \label{eq:11_golden_grad}
    \hat{H}(t)=N(\eta^D(t))y(t),
\end{equation}
which satisfies the following averaging property \footnote{The ergodic property implies that the time average of a function of the process along its trajectories exists almost surely and equals the space average:
$\lim_{T \to \infty}\frac{1}{T}\int_0^{T}a(Z_t)dt =\int_{S_Z}a(z)\mu(dz)$,  for any integrable function $a(\cdot)$, where $\mu(dz)$ is the invariant distribution of samples $Z_t$ over $S_Z$ \cite{c14}.}
\begin{equation} \label{eq:12_golden_grad}
    \frac{1}{\Pi}\int\limits_0^\Pi N(\sigma)yd\sigma=H, \quad \Pi=2\pi/\omega.
\end{equation}

Equation (\ref{eq:12_golden_grad}) was demonstrated in \cite{GKN:2012,c16,LK:2012} for the case of a quadratic nonlinear map, as in (\ref{eq:2}). Therefore, the averaged form of $H(t)$ is given by $\hat{H}_{av}(t)=(Ny)_{av}(t)=H$.

The matrix $N(\eta(t)) \in \mathbb{R}^{n\times n}$, used to obtain the Hessian estimate, is defined as:
\begin{equation} \label{eq:9}
    N_{ij}(t) = \left \{ \begin{matrix} \frac{16}{a_i^2}\Bigg(\text{sin}^2(\eta_i(t))-\frac{1}{2}\Bigg), & if \quad i=j, \\ \frac{4}{a_ia_j}\text{sin}(\eta_i(t))\text{sin}(\eta_j(t)), & if \quad i \neq j. \end{matrix} \right.
\end{equation}
%

The stochastic argument of the perturbation term is given by:
\begin{equation} \label{eq:7}
    \eta_i(t)= \omega \pi(1 + \text{sin}(W^i_{\omega t})),
\end{equation}
\noindent
representing a homogenous ergodic Markov process with a nonzero frequency $\omega>0$. As previously mentioned, the terms $W_{\omega t}^i$ indicate that distinct Wiener processes are mutually independent for each channel. By introducing the time scale $\tau = \omega t$ and applying the stochastic chain rule \cite{LK:2012}, we obtain:
\begin{equation} \label{Meq:8}
    d\eta_i =- \frac{\pi}{2}\text{sin}(W^i_{\tau})d\tau + \pi \text{cos}(W^i_{\tau})dW^i_{\tau}.
\end{equation}

Hence, from (\ref{eq:3}) and (\ref{eq:6}), the corresponding output is:
\begin{equation} \label{eq:10}
    y(t)=Q(\theta^D(t))=Q\Big(\hat{\theta}^D+S^D(t)\Big).
\end{equation}

We define the estimation error as:
\begin{equation} \label{eq:11}
    \tilde{\theta}(t):=\hat{\theta}^D(t)-\theta^\ast.
\end{equation}
Note that the error is defined using $\hat{\theta}^D$ instead of $\hat{\theta}$. With this estimation error variable, the output signal $y(t)$ can be rewritten as:
\begin{equation} \label{eq:12}
    y(t)=Q\Big( \theta^\ast+\tilde{\theta}(t)+S^D(t) \Big).
\end{equation}

To compensate for the effects of delays, we propose the following feedback law: 
\begin{align} 
    \dot{\hat{\theta}}&=U(t), \label{eq:13} \\
    \dot{U}(t)&=-cU(t)+cK\Bigg(M(\eta^D(t))y(t) \nonumber \\ 
&+N(\eta^D(t))y(t)\sum_{i=1}^n\int_{t-D_i}^tU_i(\tau)d\tau \Bigg), \label{eq:14}
\end{align}
\noindent
 for some positive constant $c\!>\!0$ and a diagonal matrix $K\!\in\! \mathrm{R}^{n \times n}$ with positive entries $K_i$, $\forall{i} \in {1,...,n}$. Since $\dot{\hat{\theta}}(t)=U^D(t)$, differentiating the error variable $\tilde{\theta}$ in (\ref{eq:11}) with respect to $t$, yields:
\begin{equation} \label{eq:15}
    \dot{\hat{\theta}}(t)=U^D(t)=\sum_{i-1}^ne_iU_i(t-D_i),
\end{equation}
\noindent
which is in the standard form of a system with input delays. Applying the reduction approach \cite{artstein} to the corresponding averaged system of (\ref{eq:15}) yields the average of the term inside the parentheses on the right-hand side of (\ref{eq:14}), corresponding to an average-based estimate $H\tilde{\theta}$ at the future time $t+D$.

\section{Stability Analysis}
The following theorem summarizes the stability and convergence properties of the closed-loop system. 
As shown in the proof, because the dynamic part is a simple integrator in 
(\ref{eq:13}) or, equivalently, in (\ref{eq:15}), the predictor feedback (\ref{eq:14}) does not require a backstepping transformation \cite{c1} to complete the stability analysis. 
The operators $\mathbb{E}\{ \cdot \}$ and $\mathbb{P} \{ \cdot \}$ denote the expected value and the probability of the signals,  respectively. 
%

\medskip
\begin{theorem} \label{maintheorem}
 Consider the closed-loop system in Figure \ref{Fig1} with multiple and distinct input delays as given in (\ref{eq:1}) and a locally quadratic nonlinear map defined by (\ref{eq:3}) and (\ref{TAChoje1}). There exists $c^{\ast} > 0$ such that, $\forall c \geq c^{\ast}$, $\exists \omega^{\ast}(c)>0$ such that, $\forall \omega \ge \omega^\ast$, the delayed closed-loop system in (\ref{eq:14}) and (\ref{eq:15}) with states $\tilde{\theta}_i(t)$, $U_i(\tau)$, $\forall{\tau} \in [t-D_i,t]$ and $\forall{i} \in \{1,\;2,\; ...\;,\;n\}$, has a unique, locally exponentially stable solution satisfying:
\begin{align} \label{eq:16}
&   \mathbb{E} \Bigg \{ \sum_{i=1}^n \Big[\tilde{\theta}_i(t)\Big]^2 + \Big [ U_i(t) \Big ] ^2 \nonumber \\ 
&   + \int\limits_{t-D_i}^t \Big [ U_i(\tau) \Big ]^2 d\tau \Bigg \}^{1/2} \leq O(1/\omega), \quad \forall{t} \rightarrow \infty.
\end{align}
In particular, 
\begin{equation} \label{eq:17}
  \lim_{(1/\omega) \to 0}\mathbb{P} \left \{ \limsup _{t \to \infty} \:|\theta(t) - \theta^{\ast}| \leq O\big ( |a|+1/\omega \big) \right \}= 1,
\end{equation}
\begin{equation} \label{eq:18}
  \lim_{(1/\omega) \to 0}\mathbb{P} \left \{ \limsup _{t \to \infty}\:|y(t) - y^{\ast}| \leq O\big ( |a|^2+1/\omega^2 \big ) \right \} = 1 ,
\end{equation}
\noindent
where $a=\big [a_1,\;a_2,\; ...\;,\; a_n \big]^T$.   
\end{theorem}
\medskip

\begin{IEEEproof} The PDE representation of the closed-loop system in (\ref{eq:14}) and (\ref{eq:15}) is given by:
\begin{equation} \label{eq:19}
    \dot{\tilde{\theta}}(t)=u(0,t),
\end{equation}
\begin{equation} \label{eq:20}
    u_t(x,t)=D^{-1}u_x(x,t), \quad x \in (0,1),
\end{equation}
\begin{equation} \label{eq:21}
    u(1,t)=U(t),
\end{equation}
\begin{align} \label{eq:22}
&    \dot{U}(t)= -cU(t)+cK \Bigg(M (\eta (t-D))y(t) \nonumber \\ 
&    + N(\eta(t-D))y(t) \int_0^1Du(x,t)dx \Bigg),
\end{align}
\noindent
where $u(x,t)=\big(u_1(x,t),u_2(x,t),\; \cdot \cdot \cdot \;,\;u_n(x,t)\big)^T \in \mathrm{R}^n$. 
Note that the solution of (\ref{eq:20}) under the condition in (\ref{eq:21}) is given by:
\begin{equation} \label{eq:23}
    u_i(x,t)=U_i(D_ix+t-D_i)
\end{equation}

\noindent
for each $i \in \{1,\;2,\; \cdot \cdot \cdot\; , n\}$. Therefore, we have:
\begin{align} \label{eq:24}
 \int_0^1D_iu_i(x,t)dx &= \int_{t-D_i}^t u_i \Bigg (\frac{\tau-t+D_i}{D_i},t\Bigg)d\tau \nonumber \\
&= \int_{t-D_i}^tU_i(\tau)d\tau.   
\end{align}

This implies that
\begin{align} \label{eq:25}
 \int_0^1Du(x,t)dx &= \sum_{i=1}^n e_i \int_0^1D_iu_i(x,t)dx \nonumber \\
&= \sum_{i=1}^n e_i \int_{t-D_i}^tU_i(\tau)d\tau.   
\end{align}

Thus, (\ref{eq:15}) can be recovered from (\ref{eq:20})-(\ref{eq:22}). The associated averaged system for (\ref{eq:19})-(\ref{eq:22}) is given by:
\begin{equation} \label{eq:26}
    \dot{\tilde{\theta}}_{av}(t)=u_{av}(0,t),
\end{equation}
\begin{equation} \label{eq:27}
    u_{av,t}(x,t)=D^{-1}u_{av,x}(x,t), \quad x \in (0,1),
\end{equation}
\begin{equation} \label{eq:28}
    u_{av}(1,t)=U_{av}(t),
\end{equation}
\begin{equation} \label{eq:29}
    \dot{U}_{av}(t)=-cU_{av}(t)+cKH\Bigg(\tilde{\theta}_{av}(t)+\int_0^1Du_{av}(x,t)dx\Bigg),
\end{equation}

\noindent
where we use the fact that the averages of $M(\eta^D(t))y(t)$ and $N(\eta^D(t))y(t)$ are given by $H\tilde{\theta}_{av}(t)$ and $H$. For simplicity in notation, let us introduce the following reduction transformation \cite{artstein}: 
\begin{equation} \label{eq:30}
    \vartheta(t):=H\Bigg(\tilde{\theta}_{av}(t)+\int_0^1Du_{av}(x,t)dx \Bigg),
\end{equation}
and define the auxiliary variable
\begin{equation} \label{eq:31}
    \tilde{U}:=U_{av}-K\vartheta.
\end{equation}

With this notation, (\ref{eq:29}) can be simplified to $\dot{U}_{av}=-c\tilde{U}$. Additionally, differentiating (\ref{eq:30}) with respect to $t$, yields:
\begin{equation} \label{eq:32}
    \dot{\vartheta}=HU_{av}(t).
\end{equation}

We establish the exponential stability of the closed-loop system by using the Lyapunov functional defined as: 
\begin{align} \label{eq:33}
&   V(t)=\vartheta(t)^TK\vartheta(t)+\frac{1}{4}\lambda_{min}(-H)\int_0^1\Big((1+x)u_{av}(x,t)^T \nonumber \\
&    \times Du_{av}(x,t)dx\Big)+\frac{1}{2}\tilde{U}(t)^T(-H)\tilde{U}(t).
\end{align}

\noindent
Recalling that $K$ and $D$ are diagonal matrices with positive entries and that $H$ is a negative definite matrix, it follows that $K$, $D$, and $-H$ are all positive definite matrices. For simplicity, the explicit dependence on time $t$ is omitted. The time derivative of $V$ is given by:
\begin{align} 
\dot{V}&=2\vartheta^TKHU_{av}+\frac{1}{2}\lambda_{min}(-H)U_{av}^TU_{av} \nonumber \\ 
&    -\frac{1}{4}\lambda_{min}(-H)u(0)^Tu(0)-\frac{1}{4}\lambda_{min}(-H) \nonumber \\
&    \times \int_0^1u_{av}(x)^Tu_{av}(x)dx+\tilde{U}^T(-H)\Big(\dot{U}_{av}-KHU_{av} \Big) \nonumber \\
& \le 2\vartheta^TKHU_{av}+\frac{1}{2}U_{av}^T(-H)U_{av}-\frac{1}{8D_{max}}\lambda_{min}(-H) \nonumber \\
&    \times \int_0^1(1+x)u_{av}(x)^TDu_{av}(x)dx \nonumber \\
&    +\tilde{U}^T(-H)\dot{U}_{av}+\tilde{U}^T(-H)K(-H)U_{av}. \label{eq:34}
\end{align}

Applying the Young's inequality to the last term yields:
\begin{align} \label{eq:35}
    \tilde{U}^T(-H)K(-H)U_{av} &\le \frac{1}{2}\tilde{U}^T(-HKHKH)\tilde{U} \nonumber \\
&    +\frac{1}{2}U_{av}^T(-H)U_{av}.
\end{align}

Then, by completing the square, we have: 
\begin{align} 
\dot{V} &\le \tilde{U}^T(-H)\tilde{U}-\vartheta^TK(-H)K\vartheta \nonumber \\ 
&    -\frac{1}{8D_{max}}\lambda_{min}(-H)\int_0^1(1+x)u_{av}(x)^TDu_{av}(x)dx \nonumber \\ 
&    +\tilde{U}^T(-H)\dot{U}_{av}+\frac{1}{2}\tilde{U}^T(-HKHKH)\tilde{U}, \nonumber \\
&   \leq\tilde{U}^T(-H)\Big(\dot{U}_{av}+c^\ast\tilde{U}\Big)-\vartheta^TK(-H)K\vartheta \nonumber \\
&    -\frac{1}{8D_{max}}\lambda_{min}(-H)\int_0^1(1+x)u_{av}(x)^TDu_{av}(x)dx, \label{eq:36}
\end{align}

\noindent
where $c^\ast:=1+\lambda_{\text{max}}(-HKHKH)/\lambda_{\text{min}}(-H)$. Thus, by setting $\dot{U}_{av}=c\tilde{U}$ for some $c>c^\ast$, it can be shown that there exists $\mu>0$ such that 
\begin{equation} \label{eq:37}
    \dot{V}=-\mu V.
\end{equation}

Finally, there exist positive constants $\alpha,\beta > 0$ such that
\begin{align} \label{eq:38}
    \alpha \Bigg(|\tilde{\theta}_{av}(t)|+\int_0^1|u_{av}(x,t)|^2dx+|\tilde{U}(t)|^2 \Bigg) \le V(t) \nonumber \\ \le \beta \Bigg (|\tilde{\theta}_{av}(t)|+\int_0^1|u_{av}(x,t)|^2dx+|\tilde{U}(t)|^2 \Bigg).
\end{align}
Thus, the averaged system in (\ref{eq:26})-(\ref{eq:29}) is exponentially stable as long as $c > c^*$. 

The procedure above eliminates the need for backstepping transformations, as used in Steps 3 and 4 of the proof of
\cite[Theorem~1]{c7}, further highlighting the simplicity of the analysis presented here. 
Essentially, while our closed-loop spectrum is finite, our non-backstepping approach does not rely on the cascade structure of a finite-spectrum system in the target variables (as would result from a backstepping transformation). Instead, it directly analyzes a feedback system without cascade decomposition.


Noting that $\int_{t-D_i}^tU_i(\tau)d\tau=\int_{-D_i}^0U_i(t+\tau)d\tau$, the closed-loop system in  (\ref{eq:14})-(\ref{eq:15}) can be rewritten as:
\begin{equation} \label{eq:39}
    \frac{d}{d\tau}\textbf{z}^{\epsilon}(\tau)=G(\textbf{z}_\tau^\epsilon)+\epsilon F(\tau,\textbf{z}_\tau^\epsilon, \eta(\tau), \epsilon),
\end{equation}
\noindent
where $\textbf{z}^\epsilon(t)=[\tilde{\theta}(t),U(t)]^T$ is the state vector, and $\epsilon:=1/\omega$. Since $\eta(\tau)$ is a homogeneously ergodic Markov process (taking values in the phase space $Y$) with invariant measure $\mu(d\eta)$ and exponential ergodicity, $\textbf{z}_t^\epsilon(\delta)=\textbf{z}^\epsilon(t+\delta)$ for $-D_n\le \delta \le 0$, and $G:\textbf{C}_2([-D_n,0]) \rightarrow{\mathbb{R}^2}$, as well as the Lipschitz function $F:\mathbb{R}_+ \times \textbf{C}_2([-D_n,0])\times Y \times [0, 1] \to \mathbb{R}^2$ with $F(\tau\,,0\,,\eta\,, \epsilon) = 0$, are continuous mappings, and $\textbf{C}_2([-D_n,0])$ denotes the class of continuous vector functions of dimension $2$ on the interval $[-D_n, 0]$, the averaging theorem by Katafygiotis and Tsarkov in \cite{c17} can be applied to conclude the exponential $p$-stability result (with $p=2$) of the initial random system, considering $\epsilon$ sufficiently small $(\epsilon \rightarrow 0)$ and thus obtain inequality (\ref{eq:16}).

From (\ref{eq:37}), the origin of the averaged closed-loop system (\ref{eq:26})-(\ref{eq:29}), incorporating the transport PDE for delay representation, is exponentially stable. Therefore, by (\ref{eq:30}) and (\ref{eq:31}), the same stability results can be concluded with respect to the norm
\begin{equation}
    \Bigg( \sum_{i=1}^n\Big[\tilde{\theta}_i^{av}(t)\Big]^2 +\int_0^{D_i}[u_i^{av}(x,t)dx]^2+[u_i^{av}(D_i,t)]^2 \Bigg) 
\end{equation}
since $H$ is non-singular.

Therefore, there exist positive constants $\bar{\alpha}$ and $\underline{\beta}$ such that all solutions satisfy $\Psi(t) \le \bar{\alpha} e^{-\underline{\beta} t}\Psi(0)$, $\forall t >0$, where $\Psi = \sum_{i=1}^n\Big[ \tilde{\theta}_i^{av}(t) \Big]^2 + \int_0^{D_i}[u_i^{av}(x,t)dx]^2 + [u_i^{av}(D_i,t)]^2$, or equivalently:
\begin{equation} \label{eq:40}
    \Psi(t) = \sum_{i=1}^n\Big[ \tilde{\theta}_i^{av}(t) \Big]^2 + \int_{t-D_i}^{t}[U_i^{av}(\tau)]^2d\tau + [U_i^{av}(t)]^2, 
\end{equation}
using (\ref{eq:23}). Thus, by the averaging theorem \cite{c17}, for a sufficiently small $\epsilon$ (or sufficiently large $\omega$), systems 
 (\ref{eq:14})-(\ref{eq:15}), or equivalently (\ref{eq:26})-(\ref{eq:29}), have a locally exponentially stable solution near the equilibrium (origin) that satisfies (\ref{eq:16}).


To prove the practical convergence to the extremum point, we define the stopping time \cite{c16}:
\begin{equation} \label{eq:53}
    \tau_{\epsilon}^{\Delta(\epsilon)}\!:=\!\inf\left \{ t \geq 0 : | \textbf{z}^{\epsilon}(t)| \!>\! {M}|\textbf{z}^{\epsilon}(0)|e^{-\lambda t} \!\!+\!\mathcal{O}(\epsilon)\right\}.
\end{equation}

\noindent
This stopping time represents the first instance when the error vector norm no longer exhibits exponential decay.
Let  
$M>0$ and $\lambda>0$ be constants, and let $T(\epsilon):(0,1) \rightarrow \mathbb{N}$ be a continuous function. Then, analogously to \cite{c16}, the norm of the error vector $|\textbf{z}^{\epsilon}(t)|$ converges both \emph{almost surely} (a.s.) and \emph{in probability} to a value below a residue $\Delta(\epsilon)=\mathcal{O}(\epsilon)$:
\begin{equation} \label{eq:sm1}
    \lim_{\epsilon \to 0} \inf \{ t \geq 0 : |\textbf{z}^{\epsilon}(t)| > M |\textbf{z}^{\epsilon}(0)|e^{-\lambda t} + \Delta \} = \infty , \textrm{a.s.},
\end{equation}

\begin{equation} \label{eq:sm2}
   \lim_{\epsilon \to 0} \mathbb{P}\{|\textbf{z}^{\epsilon}(t)| \leq M|\textbf{z}^{\epsilon}(0)|e^{-\lambda t} + \Delta\,, \forall{t} \in [0,T(\epsilon)]\}=1,
\end{equation}

\noindent
with $\lim_{\epsilon \to 0}\textrm{T}(\epsilon) = \infty$. From (\ref{eq:sm1}), it is evident that $\tau_{\epsilon}^{\Delta (\epsilon)}$ surely tends to infinity as $\epsilon$ approaches zero. Similarly, in (\ref{eq:sm2}), the deterministic function $T(\epsilon)$ also tends to infinity as $\epsilon$ approaches zero. This implies that exponential convergence holds over an arbitrarily long time interval, ensuring that each component of the error vector converges below the residue $\Delta(\epsilon)=\mathcal{O}(\epsilon)$, specifically $\Tilde{\theta}(t)$. Consequently, we have $\lim_{\epsilon \to 0}\mathbb{P} \Big \{ \textrm{lim sup}_{t \to \infty}$ $|\tilde{\theta}(t)|\Big\} \leq \mathcal{O}(\epsilon) \Big \}= 1$. From equations (\ref{eq:11}) and (\ref{eq:6}), it follows that:
\begin{equation} \label{eq:sm11}
  \theta^D(t)- \theta^{\ast} =  \tilde{\theta}(t)+S(\eta^D(t)). 
\end{equation}

Since the first term on the right-hand side of (\ref{eq:sm11}) is ultimately of order  $\mathcal{O}(\epsilon)$, and the second term is of order $\mathcal{O}(|a|)$, we obtain equation (\ref{eq:17}). Finally, by combining (\ref{eq:3}) and (\ref{eq:17}), we derive (\ref{eq:18}). 
\end{IEEEproof} 

A deterministic version of the results presented in Theorem~\ref{maintheorem} can be found in the companion conference paper \cite{c21}.  


\section{Simulations}
To evaluate the effectiveness of delay compensation in stochastic multivariable ESC, we consider a static quadratic map defined as:
\begin{equation} \label{eq:100}
    Q(\theta)=5+\frac{1}{2}\Big(2(\theta_1)^2+4(\theta_2-1)^2+4\theta_1(\theta_2-1) \Big), 
\end{equation}

\noindent
with input delays $D_1=50$ and $D_2=100$, optimizers located at $\theta^\ast=[0\,\,1]^T$, and an extremum point $y^\ast=5$. The Hessian of the map is given by 
\begin{equation} \label{eq:Hessiana}
    H=-\begin{pmatrix} 
        2 &\; 2 \\ 
        2 &\; 4 \\
        \end{pmatrix}.
\end{equation}

In this case ($n=2$), the predictor equation simplifies to:
\begin{align} \label{eq:ContPred}
&    \dot{U}(t)= -cU(t)+c \begin{bmatrix} 
  K_1 &  0 \\ 
   0  & K_2 \\
  \end{bmatrix} \Big( M(\eta^D(t))y(t)+N(\eta^D(t))y(t) \nonumber \\
&  \times \Bigg( \begin{bmatrix} 
  1 \\
  0  \\
  \end{bmatrix} \int_{t-D_i}^t U_1(\tau)d\tau + \begin{bmatrix}
  0 \\
  1  \\
  \end{bmatrix} \int_{t-D_2}^t U_2(\tau) d\tau \Bigg) \Big).
\end{align}

Given that $M(\eta(t))$ and $N(\eta(t))$ are defined as: 
\begin{equation} \label{eq:M}
M(\eta^D(t))=\begin{bmatrix}
    M_1(\eta(t-D)) \\
    M_2(\eta(t-D)) 
\end{bmatrix} = \begin{bmatrix}
\frac{2}{a_1}\textrm{sin}(\eta_1(t-D_1)) \\
\frac{2}{a_2}\textrm{sin}(\eta_2(t-D_2))
\end{bmatrix},
\end{equation}
and 
\begin{align} \label{eq:N}
& N(\eta^D(t))=\begin{bmatrix}
N_{11}(\eta^D(t)) & N_{12}(\eta^D(t)) \\
N_{21}(\eta^D(t)) & N_{22}(\eta^D(t))  \\
\end{bmatrix} = \nonumber \\
& \begin{small}
\begin{bmatrix}
\hspace{-2.0 cm}-\frac{8}{a^2}\text{cos}^2(2\eta_1(t-D_1)) &\hspace{-1.8 cm} \frac{4}{a_1a_2}\text{sin}(\eta_1(t-D_1))\text{sin}(\eta_2(t-D_2)) \\
\frac{4}{a_1a_2}\text{sin}(\eta_1(t-D_1))\text{sin}(\eta_2(t-D_2)) & -\frac{8}{a^2}\text{cos}^2(2\eta_2(t-D_2))  \\
\end{bmatrix},
\end{small}
\end{align}
the predictor equation (\ref{eq:ContPred}) can be separated into expressions for $\dot{U}_1(t)$ and $\dot{U}_2(t)$ as follows:
\begin{align} \label{eq:U1}
&    \dot{U}_1(t)\!=\!-cU_1(t)\!+\!cK_1\Big( M_1(\eta(t\!-\!D))y(t)\!+\!y(t) \Big( N_{11}(\eta^D(t)) \nonumber \\ 
&     \times \int_{t-D_1}^tU_1(\tau)d\tau \!+\!N_{12}(\eta^D(t))\int_{t-D_2}^tU_2(\tau)d\tau \Big) \Big),
\end{align}

\begin{align} \label{eq:U2}
&    \dot{U}_2(t)\!=\!-cU_2(t)\!+\!cK_2\Big( M_2(\eta(t\!-\!D))y(t)\!+\!y(t) \Big( N_{21}(\eta^D(t)) \nonumber \\ 
&    \times \int_{t-D_1}^tU_1(\tau)d\tau \!+\!N_{22}(\eta^D(t))\int_{t-D_2}^tU_2(\tau)d\tau \Big) \Big).
\end{align}
For comparison, note that the expressions for $U_1(t)$ and $U_2(t)$ in (\ref{eq:U1}) and (\ref{eq:U2}) are significantly simpler than the deterministic control laws proposed in \cite[Equations (49) and (160)-(161)]{c7}, specifically for the case $n=2$. 

The additive dither signal is defined as: 
\begin{equation}
S(t) = 
	\begin{bmatrix}
		a_1 \sin(\eta_1(t)) &
		a_2 \sin(\eta_2(t))
	\end{bmatrix}^T.
\end{equation}    


Numerical simulations of the predictor in (\ref{eq:ContPred}) were conducted using the following parameter values: $a_1=a_2=0.22$, $c=20$, $\omega=5$, $K=0.005\; \text{I}$, where $\text{I} \in \mathbb{R}^{2\times2}$ is the  identity matrix, and  $\hat{\theta}(0)=(1,0)$.

\begin{figure}[htb!]
    \centering
    \includegraphics[scale=0.5]{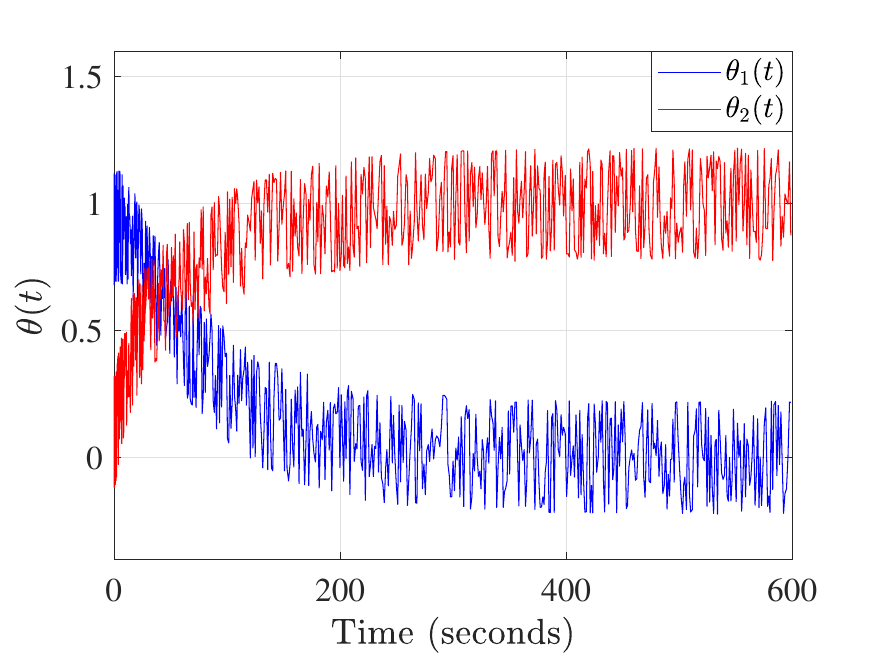}
    \caption{System input $\theta(t)$ for the multivariable gradient-based stochastic ESC in the absence of delays.}
    \label{fig:2}
\end{figure}

\begin{figure}[htb!]
    \centering
    \includegraphics[scale=0.5]{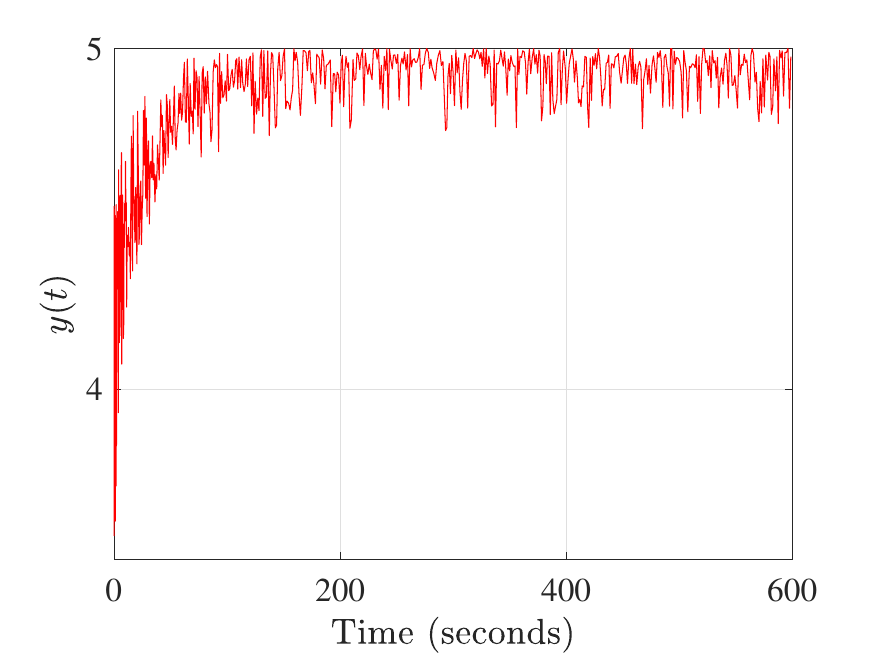}
    \caption{System output \( y(t) \) for the multivariable gradient-based stochastic ESC in the absence of delays.}
    \label{fig:3}
\end{figure}

Figure \ref{fig:2} illustrates 
the convergence of the system input to a neighborhood around the optimizer $\theta^\ast=[0 \,\,1]^T$ 
when delays are not considered, and the classical gradient-based stochastic approach \cite{LK:2012, c16} is applied without a prediction scheme.
Additionally, Figure \ref{fig:3} shows the output of the nonlinear map converging to a small vicinity of the extremum point $y^\ast=5$. 

\begin{figure}[htb!]
    \centering
    \includegraphics[scale=0.5]{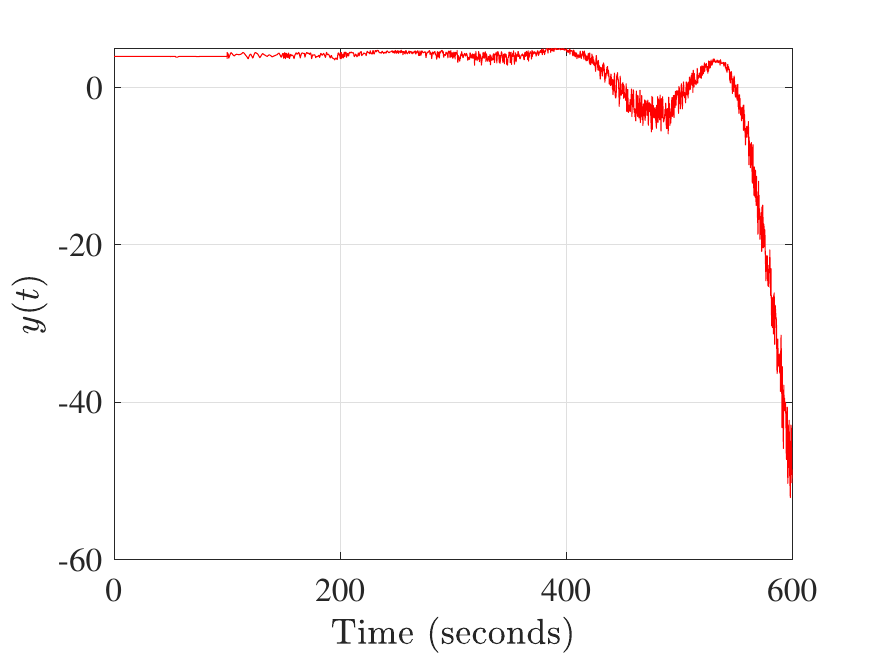}
    \caption{System output $y(t)$ for the multivariable gradient-based stochastic ESC with delays, without predictor feedback.}
    \label{fig:5}
\end{figure}

In Figure \ref{fig:5}, 
it can be observed that the classical stochastic extremum-seeking approach lacks robustness against the introduced delays, resulting in closed-loop system instability, as anticipated from the literature \cite{c7}. 
This highlights the importance of incorporating the proposed predictor-feedback scheme within the stochastic ESC framework. 

\begin{figure}[htb!]
    \centering
    \includegraphics[scale=0.5]{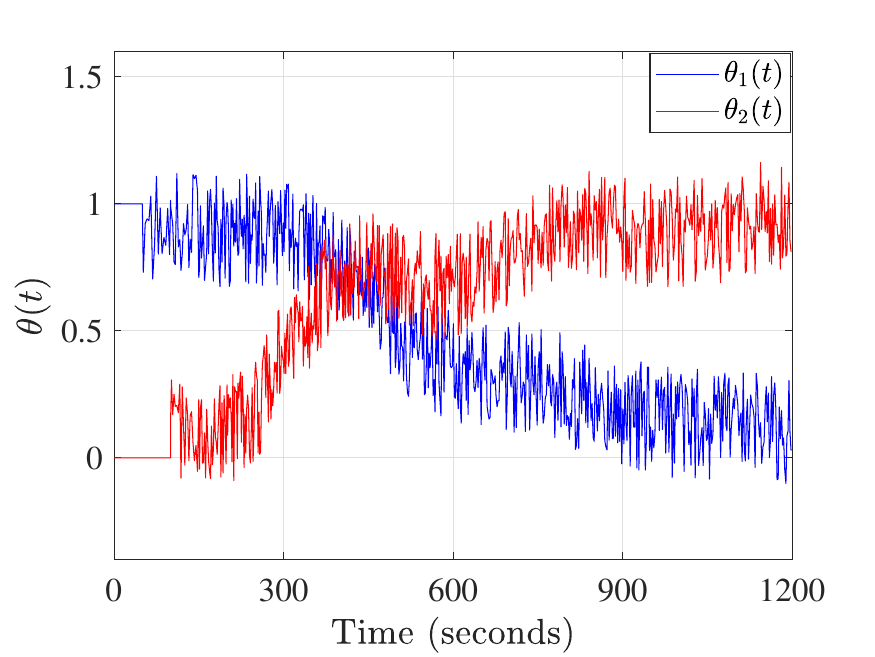}
    \caption{System input $\theta(t)$ for the multivariable gradient-based stochastic ESC with delays and predictor feedback.}
    \label{fig:6}
\end{figure}

\begin{figure}[htb!]
    \centering
    \includegraphics[scale=0.5]{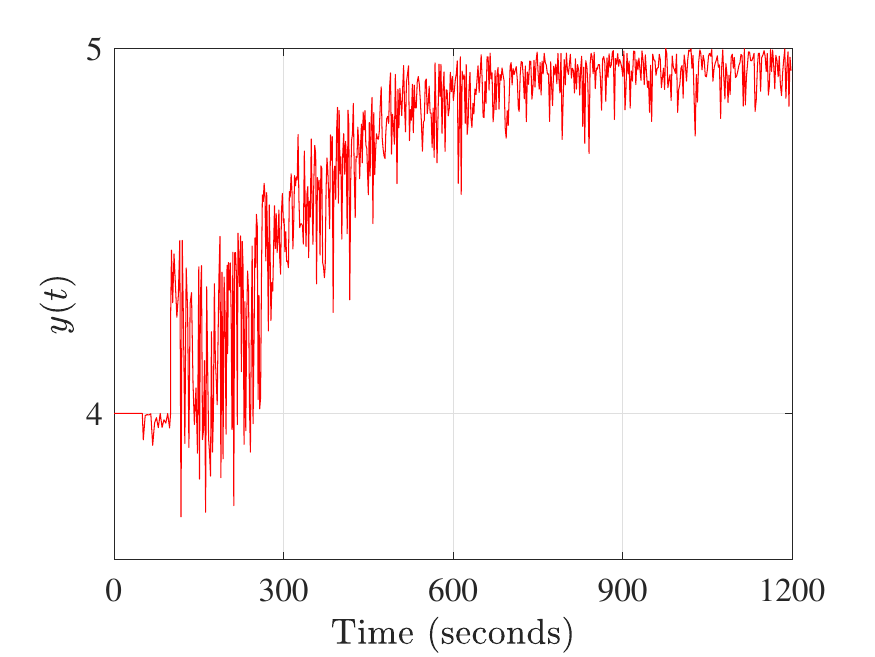}
    \caption{System output $y(t)$ for the multivariable gradient-based stochastic ESC with delays and predictor feedback.}
    \label{fig:7}
\end{figure}

\begin{figure}[htb!]
    \centering
    \includegraphics[scale=0.5]{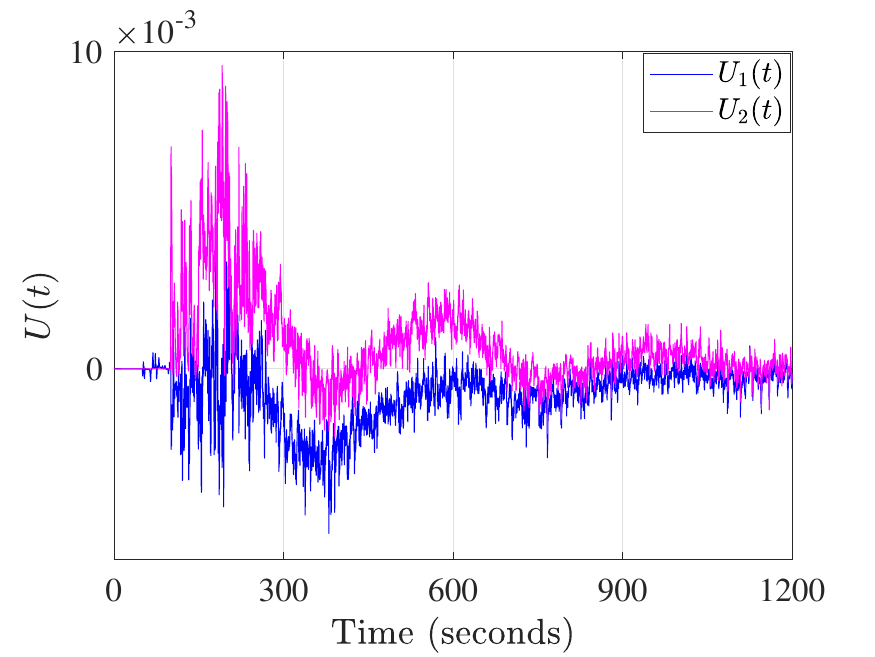}
    \caption{Control signal $U(t)$ for the multivariable gradient-based stochastic ESC with delays and predictor feedback.}
    \label{fig:8}
\end{figure}

\begin{figure}[htb!]
    \centering
    \includegraphics[scale=0.5]{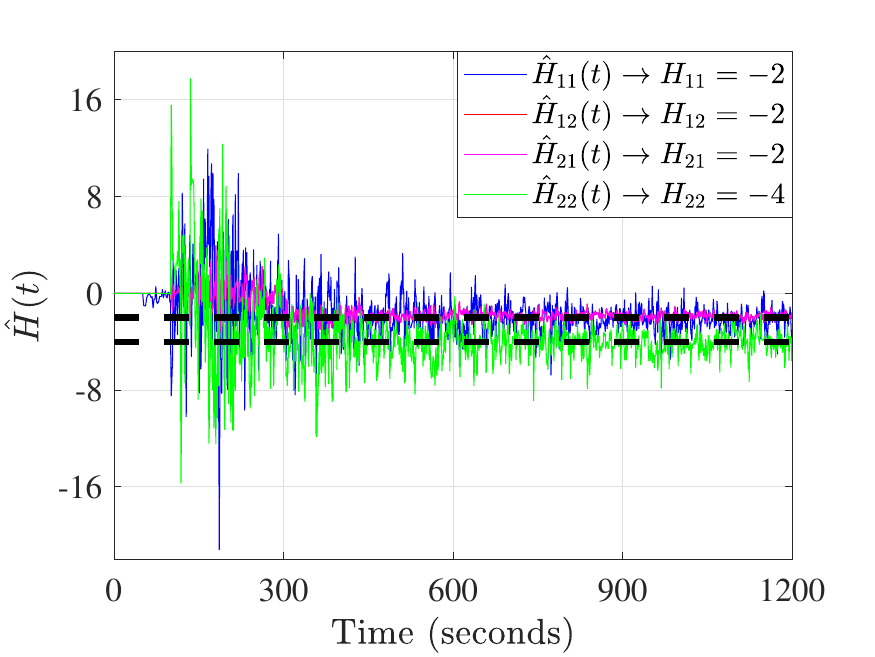}
    \caption{Convergence of the elements of $\hat{H}(t)$ to the Hessian matrix $H$. The first black dashed line corresponds to $H_{11}=H_{12}=H_{21}=-2$, and the second black dashed line corresponds to $H_{22}=-4$.}
    \label{fig:9}
\end{figure}

Figure \ref{fig:6} confirms the preservation of closed-loop system stability and the convergence of the plant input to a neighborhood around the optimal value $\theta^{\ast}=[0 \,\, 1]^T$ when the proposed predictor-feedback scheme is applied. 
In Figure \ref{fig:7}, the predictor is shown to effectively compensate for the inserted delays, allowing the plant output to converge toward the desired value $y^{\ast}=5$. Additionally, Figure \ref{fig:8}, illustrates the attenuation of the control signal $U(t)$ as the output of the map approaches the neighborhood of its extremum value $y^\ast=5$.
Finally, Figure \ref{fig:9}, shows that the estimate $\hat{H}(t)$ of the Hessian matrix in (\ref{eq:11_golden_grad}) converges to the true values of the unknown elements of $H$ in (\ref{eq:Hessiana}).

%
%

\section{Conclusion}
In this paper, we proposed a new multivariable gradient-based stochastic extremum-seeking control scheme for real-time optimization of multi-input static maps with distinct actuator delays. This is a remarkable advance since our previous results \cite{JCAES,CODIT} involving multiparameter stochastic extremum-seeking control was limited to equal delays in the input channels. The introduced control law compensates for multiple and distinct delays by employing predictor feedback with averaging-based Hessian estimation, along with appropriately tuned, time-shifted stochastic dither signals. 
This predictor-based feedback control law is notably simpler than previous deterministic gradient-based approaches in the literature \cite{c7}. Furthermore, unlike other methods relying on sequential predictors \cite{A16} or sampled-data designs \cite{Emilia_arxiv}, our approach imposes no restrictions on delay duration to sustain convergence rates. Instead, it achieves reliable cross-coupling predictions across channels with arbitrarily long (and distinct) input delays, maintaining convergence speed.

A key advancement in this framework is the use of Artstein's classical reduction approach \cite{artstein} rather than the infinite-dimensional backstepping methodology \cite{c1} typically used in delay-compensated extremum-seeking studies. This innovation enhances the practicality and implementation simplicity of the proposed scheme, offering a streamlined and effective solution to the challenges of multivariable stochastic extremum-seeking in systems with delay.

In future work, it would be worthwhile to explore the handling of distributed delays \cite{c223} or other classes of PDEs, as discussed in \cite{white_book}. For instance, the authors of \cite{NewHigh} presented a Newton-based stochastic ESC algorithm that focuses on maximizing higher derivatives of unknown maps under output delays, introducing stochastic perturbations to accommodate arbitrarily long delays in dynamic single-input maps. Extending this approach to multi-input maps with distinct input delays would be an intriguing direction for further investigation, potentially broadening the applicability and robustness of stochastic extremum-seeking control in complex delayed systems.


\bibliographystyle{IEEEtranS}
\bibliography{IEEEabrv,references}

\end{document}